\documentclass{article}

\usepackage[all]{xy}
\usepackage{amssymb}

\newcommand{\cat}{\mathbf}

\newcommand{\C}{\ensuremath{\mathbf{C}}}

\newcommand{\R}{\ensuremath{\mathbf{R}}}

\newcommand{\rad}{\ensuremath{\mathrm{rad}}}

\newcommand{\qed}{\hfill\strut\nobreak\hfill\nobreak$\square$}

\newcommand{\ie}{\textit{i.e.}}

\newtheorem{theorem}{Theorem}
\newtheorem{proposition}[theorem]{Proposition}

\title{Variations on a theme of Gel'fand and Na\u\i{}mark}
\author{Miguel Carri\'on \'Alvarez\\
        Department of Mathematics\\
        University of California\\
        Riverside, CA 92521, USA\\
        \texttt{miguel@math.ucr.edu}}

\begin{document}

\maketitle

\begin{abstract}

$C^*$-algebras are widely used in mathematical physics to represent
the observables of physical systems, and are sometimes taken as the
starting point for rigorous formulations of quantum mechanics and
classical statistical mechanics. Nevertheless, in many cases the
na\"\i{}ve choice of an algebra of observables does not admit a
$C^*$-algebra structure, and some massaging is necessary. In this
paper we investigate what properties of~$C^*$-algebras carry over to
more general algebras and what modifications of the Gel'fand theory of
normed algebras are necessary. We use category theory as a guide and,
by replacing the ordinary definition of the Gel'fand spectrum with a
manifestly functorial definition, we succeed in generalizing the
Gel'fand--Na\u\i{}mark theorem to locally convex $*$-algebras.  We
also recall a little-known but potentially very useful generalization
of the Stone--Weierstrass theorem to completely regular, Hausdorff
spaces.

AMS Mathematics Subject Classification (2000): 46M99 (primary) 47L90
(secondary).

\end{abstract}

\section{Introduction}

In rigorous formulations of quantum theory, non-commutative
$C^*$-algebras are used extensively to represent the observables of
physical systems~\cite{haag96,emch72,BSZ,wald94}. This approach was
pioneered by Segal~\cite{segal47}, who also advocated formalizing
probability theory in terms of commutative algebras of bounded random
variables~\cite{MR16:149d}; this leads naturally to a $C^*$-algebraic
formulation of classical statistical mechanics~\cite{emch72}. Once the
$C^*$-algebra of observables is specified, the formal development of
either theory is well-understood. In the case of non-commutative
$C^*$-algebras, the Gel'fand--Na\u\i{}mark--Segal construction
produces representations of the $C^*$-algebra as bounded operators on
a Hilbert-space, and the problem is to identify and analyze the
physically relevant ones. In the commutative case, the Gel'fand
transform provides a geometric interpretation of the algebra of
observables as continuous bounded functions on a compact space, and
the Gel'fand--Na\u\i{}mark--Segal construction gives this algebra of
random variables a Hilbert-space structure with the covariance as
inner product.

In physics, the hardest problem is often finding a suitable
$C^*$-algebra of observables in the first place, when the only input
is a geometrical or operational description of a physical system.  As mentioned before, the Gel'fand theory of
normed algebras provides a natural interpretation of Abelian
$C^*$-algebras as algebras of continuous functions on a compact
Hausdorff space, the Gel'fand spectrum, which in applications to
classical mechanics would be the phase space of the system under
consideration. It should be obvious, however, that compact phase
spaces are very rare as, even when the configuration space is bounded,
momentum is usually unbounded. Even worse, in the case of field
theories or the mechanics of continuous media, the phase space is an
infinite-dimensional manifold and is not even locally compact! Also,
in many instances, such as when using Poisson brackets, one is
interested in algebras of smooth functions which, while tending
to be metrizable, do not admit a norm and so cannot be
$C^*$-algebras.

Given an algebra of classical observables which is not a
$C^*$-algebra, two related approaches can be taken. The first is to
study the manipulation required to turn the given algebra into a
$C^*$-algebra. The process may involve loss of information (as when
taking equivalence classes), new information (as when extending or
completing the algebra), or arbitrary choices (such as a choice of
complex structure on a real space), and one should pay attention to
the physical interpretation of these manipulations. The second
approach, which we develop in this paper, is trying to extend as much
as possible of the theory of $C^*$-algebras to more general algebras,
possibly changing the key definitions in the theory so they apply more
generally. For the mathematical side of this exploration we use
category theory, which provides notation and concepts tailored to
asking and answering questions about naturality of mathematical
operations. As for the physical interpretation, whether a manipulation
is unphysical can only be answered in each particular instance, but
hopefully mathematically sensible manipulations will turn out to be
physically sensible.

In the Gel'fand theory of normed algebras~\cite{rickart60}, the
Gel'fand spectrum of maximal ideals of an algebra plays a central
r\^ole. As we have indicated, it is given a natural compact, Hausdorff
topology, and the elements of the algebra can be naturally interpreted
as continuous complex functions on it. One of the key results in the
theory of $C^*$-algebras is the Gel'fand--Na\u\i{}mark theorem, which
states that every Abelian $C^*$-algebra is isometrically
$*$-isomorphic to the $*$-algebra of bounded continuous functions on
its spectrum, with the normed topology of uniform convergence. The
content of the Gel'fand--Na\u\i{}mark theorem really is that the
quotient of an Abelian $C^*$-algebra by a maximal ideal is a
continuous $*$-homomorphism into the complex numbers. This is a
striking connection between algebra and analysis, but the functorial
properties of the Gel'fand transform depend not on the fact that the
spectrum consists of maximal ideals, but on the fact that it is a
hom-set in the category of $C^*$-algebras.  Accordingly, in more
general settings than $C^*$-algebras it is more productive to simply
restrict one's attention to continuous $*$-homomorphisms into the
complex numbers than to study maximal ideals. Large parts of the
Gel'fand--Na\u\i{}mark theory, including the
Gel'fand--Na\u\i{}mark--Segal construction, then carry over.

Our main results, theorems~\ref{thm:spectrum} and~\ref{thm:transform},
imply that, under rather general hypothesis, given a $*$-algebra~$A$
one can find a $*$-homomorphism injecting it as a dense $*$-subalgebra
of continuous complex functions on a Tychonoff (completely regular,
Hausdorff) space, with the compact-open topology. This
$*$-homomorphism always exists, but it may have a nontrivial kernel.

The contents of this paper are summarized in the following diagram of
functors, each cell of which roughly corresponds to one section:
$$
\xymatrix{\cat{AbAlg}\ar@<.5ex>[r]^{\mathrm{ F}}\ar[d]_{\omega(-,-^*)}&\cat{Ab{}^*Alb}\ar@<.5ex>[l]^{\mathrm{U}}\ar[d]^{\omega(-,-^*)}\\
          \cat{LCAbAlg}\ar@<.5ex>[r]^{\mathrm{ F}}\ar[d]_{\Delta_\cdot}&\cat{LCAb{}^*Alg}\ar@<.5ex>[l]^{\mathrm{U}}\ar[dl]_{\Delta_\cdot}\ar[d]^{A\mapsto\bar{\hat
          A}}\\
          \cat{Tych}\ar@<.5
          ex>[r]^{C(-)}&\cat{AbLC^*Alg}\ar@<.5
          ex>[l]^{\Delta_\cdot}\\}
$$

In section~\ref{sec:AbAlg} we define the categories~$\cat{AbAlg}$
and~$\cat{Ab{}^*Alb}$ of unital Abelian algebras and $*$-algebras, and
study the adjoint pair of functors, ``underlying'' and ``free'', going
between them.

Section~\ref{sec:WkTop} deals with the square cell at the top of the
above diagram of functors. We use the weak topology to make every
algebra and~$*$-algebra locally convex, hence the
names~$\cat{LCAbAlg}$ and~$\cat{LCAb{}^*Alg}$, in such a way that the
underlying and free functors commute with the operation of adding the
weak topology.

The Gel'fand spectrum is defined as a hom-set in
section~\ref{sec:spectrum}, associated to the triangular cell to the
left of the diagram of functors. This definition entails that the
spectrum of a~$*$-algebra is, in general, strictly contained in the
spectrum of its underlying algebra. It is then shown that the Gel'fand
spectrum is a Tychonoff space, the next best thing after compact
Hausdorff spaces, and that it is a weak${}^*$-closed subset of the
topological dual of the algebra. Because of the functorial definition
of the spectrum, the Gel'fand transform automatically becomes
a~$*$-algebra $*$-homomorphism when applied to a~$*$-algebra, in a
setting when the lack of a norm makes the usual techniques applied
to~$C^*$-algebras break down. We also discuss the interpretation of
the statement ``every $*$-homomorphism is a homomorphism'' as a
natural transformation.

The triangular cell on the bottom-right of the diagram is discussed in
section~\ref{sec:trans}. We study the difference between the usual
definitions of the Gel'fand transform and ours, and use a
generalization of the Stone--Weierstrass theorem from the case of a
compact, Hausdorff space (the spectrum of a $C^*$-algebra) to the case
of a Tychonoff space (the spectrum of any algebra), to show that the
image of the Gel'fand transform is dense in the continuous functions
on the spectrum. The generalization of the Stone--Weierstrass theorem
involves replacing uniform convergence by the compact-open topology
(uniform convergence on compact sets), which is not a surprise since
this is a topology that is extensively used in complex analysis. For
lack of a better name, we call the space of complex continuous
functions on a Tychonoff spaces an ``Abelian $LC^*$-algebra'' ($LC$
for locally convex, or for ``locally~$C^*$'').

Finally, in section~\ref{sec:GNS} we study the notion of a state and
apply the Gel'fand--Na\u\i{}mark--Segal construction to
the~$*$-algebra of complex functions on a Tychonoff space. The states
are realized as compactly-supported Borel probability measures, which
is related to the fact that the restriction of the algebra of
continuous functions to a compact set is a~$C^*$-algebra. This
illustrates the sense in which we are dealing with ``locally $C^*$''
algebras.

While lacking an intrinsic, algebraic characterization
of~$LC^*$-algebras (such as is available for $C^*$-algebras), our
discussion shows that there is life outside the world of
$C^*$-algebras in the sense that the basic operations that
mathematical physicists need to perform on algebras of observables can
be carried out for~$LC^*$-algebras. Also, it illustrates how category
theory can be a powerful guide to find the right definitions making it
possible to extend impressive results like the Gel'fand--Na\u\i{}mark
theory to situations where few, if any, of the specific techniques
used in the original proofs are available.

\section{Abelian algebras and $*$-algebras}\label{sec:AbAlg}

For the purposes of this paper, an algebra will be a complex vector
space~$A$ with an associative, bilinear multiplication and a
unit~$1_A$. An algebra is Abelian if~$ab=ba$ for all~$a,b\in A$. A
linear map~$\phi\colon A\to B$ between algebras is an algebra
homomorphism if, and only if, $\phi(aa')=\phi(a)\phi(a')$ and
$\phi(1_A)=1_B$. Note that we are assuming that all algebras are
unital, and that all algebra homomorphisms map units to units. This is
partly because the category of unital algebras with unit-preserving
homomorphisms is relatively nice among the possible categories of
algebras. We denote this category by~$\cat{Alg}$, and the category of
Abelian algebras with algebra homomorphisms by~$\cat{AbAlg}$.

An algebra~$A$ is a $*$-algebra if it has an involutive anti-linear
anti-homomorphism~$*\colon A\to A$. What this means is that
$a^{**}=a$, that $(z1_A)^*=\bar z1_A$, and that
$(ab)^*=b^*a^*$. Although in general the term ``involution'' refers
just to~$*$ being its own inverse, in this context a $*$-algebra is
usually called ``an algebra with an involution'', and the
operation~$*\colon a\mapsto a^*$ is called ``the involution''. The
complex numbers is naturally a $*$-algebra, with involution given by
complex conjugation,~$z^*\colon =\bar z$.

An algebra homomorphism $f\colon A\to B$ between $*$-algebras is an
algebra $*$-homomorphism if, and only if, $f(a^*)=f(a)^*$ for all
$a\in A$. Note that our definition of the involution includes the
requirement that the ``unit map''~$e_A\colon\C\to A$, such
that~$e_A(z)=z1_A$, be a~$*$-homomorphism. There is a category of
$*$-algebras and $*$-homomorphisms, which we denote~$\cat{Alg}^*$, and
a category of Abelian $*$-algebras with $*$-homo\-morph\-isms
denoted~$\cat{Ab{}^*Alb}$.

It is clear that every (Abelian) $*$-algebra is an (Abelian) algebra,
and that every $*$-homomorphism is a homomorphism. Hence, the process
of considering a $*$-algebra as an algebra is a forgetful
functor~$U\colon\cat{Alg}^*\to\cat{Alg}$ restricting to a functor
between the categories of Abelian algebras~$
U\colon\cat{Ab{}^*Alb}\to\cat{AbAlg}$. If~$A$ is an (Abelian)
$*$-algebra,~$ U(A)$ is called its underlying (Abelian) algebra. As it
happens generally in algebra, these forgetful functors have
left-adjoint functors. In the Abelian case~$
F\colon\cat{AbAlg}\to\cat{Ab{}^*Alb}$ is such that there is a natural
isomorphism
$$
\cat{AbAlg}\bigl(A, U(B)\bigr)\simeq
\cat{Ab{}^*Alb}\bigl( F(A),B\bigr)
$$
for all Abelian algebras~$A$ and Abelian $*$-algebras~$B$
(following~\cite{lane98}, we denote by~$\cat{Xmpl}(A,B)$ the set of
morphisms~$f\colon A\to B$ in the category~$\cat{Xmpl}$). The
functor~$ F$ is said to be the left adjoint of~$ U$, and its
interpretation is that~$ F(A)$ is the free Abelian $*$-algebra
generated by~$A$.

The existence of the functor~$ F$ is equivalent to the following
universal property: for every Abelian algebra~$A$ there exists an
Abelian $*$-algebra~$ F(A)$ and such that, for every Abelian
$*$-algebra~$B$, and for every algebra homomorphism~$f\colon A\to
U(B)$, there exists a unique $*$-homomorphism~$f'\colon F(A)\to B$
such that the following diagram commutes:
$$
\xymatrix{A\ar[r]^{\iota_A}\ar[dr]^f& U\bigl(F(A)\bigr)\ar[d]^{U(f')}\\
          &U(B)\\}
$$
where~$\iota_\cdot$ is ``the unit of the adjunction''. Uniqueness
of~$F(A)$ up to isomorphism follows from abstract
nonsense~\cite{lane98}.

Existence of~$ F(A)$ is proved by construction. In the Abelian case,
it suffices to add to~$A$ a new element~$a^*$ for each $a\in A$ and
all the products of the form~$a^*b$, and to consider the collection of
all linear combinations of those three kinds of elements, subject only
to the relations necessary to enforce that~$a\mapsto a^*$ is an
involution on the new algebra. Note that, although the ``underlying
algebra'' functor is the same for Abelian and non-Abelian algebras,
the ``free algebra'' functor is very different from the ``free Abelian
algebra'' functor. In particular, the free $*$-algebra generated by an
Abelian algebra is non-Abelian, and much larger than the free Abelian
$*$-algebra generated on it, because n the free Abelian $*$-algebra,
$a^*b=ba^*$ but not so in the free $*$-algebra. Finally, note that
neither~$ U F$ nor~$ F U$ are the identity. We summarize the situation
thus:
$$
\xymatrix{\cat{AbAlg}\ar@<.5ex>[r]^{\mathrm{
F}}&\cat{Ab{}^*Alb}\ar@<.5ex>[l]^{\mathrm{ U}}\\}
$$

\subsection{Examples}

The simplest example of this involves algebras of polynomials with
complex coefficients. 
\begin{description}
\item[{$\C[x]$}]Let us start by considering the $*$-algebra of complex
polynomials on one self-adjoint variable~$x$ (satisfying~$x^*=x$), so
the involution maps a polynomial~$a_0+a_1x+\cdots+a_nx^n$
to~$\overline{a_0}+\overline{a_1}x+\cdots+\overline{a_n}x^n$. We
denote this $*$-algebra by~$A=\C[x]$.

In application to classical mechanics, this is the polynomials on one
real configuration variable.

\item[{$\C[z]$}]The underlying algebra~$B= U(A)$ is the same
algebra of polynomials, except that it is ``forgotten'' that one can
apply the involution to them. So, the two
polynomials~$a_0+a_1z+\cdots+a_nz^n$ and
$\overline{a_0}+\overline{a_1}z+\cdots+\overline{a_n}z^n$ are both
elements of~$B$, but now they are not related by any operation on~$B$
unless all the~$a_i$ are real and the polynomials are actually the
same. We denote this algebra by~$B=\C[z]$, which is actually the
``free Abelian algebra on one generator''.

In classical mechanics, this is the algebra of polynomials on one
complex phase variable~$z=q+ip$.

\item[${\C[z,z^*]}$]We now consider~$C=F(B)$, the free~$*$-algebra
on~$B$ or the ``free~$*$-algebra on one generator''. To~$z$ we must
add a distinct adjoint~$z^*$, and then build the free Abelian algebra
generated by the two. A typical polynomial in this algebra is
on~$a_{00}+a_{10}z+a_{01}z^*+a_{20}z^2+a_{11}zz^*+a_{02}(z^*)^2+\cdots$,
and the effect of the involution on it is
now~$\overline{a_{00}}+\overline{a_{10}}z^*+\overline{a_{01}}z+\overline{a_{20}}(z^*)^2+\overline{a_{11}}zz^*+\overline{a_{02}}z^2+\cdots$. We
denote this $*$-algebra by~$C=\C[z,z^*]$.

In classical mechanics this is the algebra of polynomials on phase
space, since we can interpret~$z=q+ip$ and $z^*=q-ip$ with $p=p^*$ and
$q=q^*$, and indeed in that case 
$$
\C[z,z^*]\simeq\C[q,p].
$$
where~$q$ and~$p$ are self-adjoint generators. Accordingly, in
physical applications this algebra could also be associated to a
two-dimensional real configuration space.

\item[{$\C[z,w]$}]The underlying algebra of~$\C[z,z^*]$
is~$D=\C[z,w]$, where we still have two generators but we forget that
they are related by the involution (or, equivalently, we forget that
the two generators are self-adjoint). This is algebra of complex
polynomials on two variables.

In classical mechanics, this would be the algebra of polynomials on a
phase space of two degrees of freedom, with~$z=q_1+ip_1$
and~$w=q_2+ip_2$.

\end{description}
The situation is summarized by the following diagram:
$$
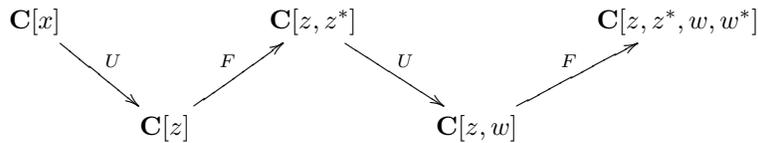
\xymatrix{\C[x]\ar[rd]^{ U}&&\C[z,z^*]\ar[rd]^{ U}&&\C[z,z^*,w,w^*]\\
          &\C[z]\ar[ru]^{ F}&&\C[z,w]\ar[ru]^{ F}\\}
$$
It is apparent that, in this example,~$F$ doubles the complex
dimension of the algebra as a vector space, while~$ U$ leaves it
unchanged.

Another interesting series of examples, this time related to quantum
field theory, is that of algebras where the generators form a Hilbert
space. These are associated to the Fock space representation of
systems with variable numbers of particles, such as are used in
particle physics, quantum optics or solid-state physics. As algebras,
they are the algebras of creation operators, which are Abelian
subalgebras of the full-blown algebras of observables on Fock space.

Let~$H$ be a complex Hilbert space. If~$H$ has an anti-unitary
involution~$*$, $H$ decomposes as~$H\simeq H^\sharp\oplus_\R
iH^\sharp$, where~$H^\sharp$ is the real eigenspace of vectors such
that~$a^*=a$. We denote the polynomials on~$H$ by~$\C[H^\sharp]$
if~$H$ has an involution, and~$\C[H]$ otherwise. The analogue of the
preceding diagram is
$$
\xymatrix{\C[H^\sharp]\ar[rd]^{ U}&&\C[H\oplus
          H^*]\simeq\C[H_1^\sharp\oplus H_2^\sharp]\ar[rd]^{ U}\\
          &\C[H]\ar[ru]^{ F}&&\C[H_1\oplus H_2]\\}
$$
and the physical interpretation of each of the algebras is as follows.
\begin{description}
\item[${\C[H^\sharp]}$]The space~$H^\sharp$ is a real Hilbert space of
states of a truly neutral particle (which is its own antiparticle),
and its complexification~$H$ is the complex vector space of all
quantum states. The algebra~$\C[H^\sharp]$ is a dense subspace of the
corresponding Fock space or, alternatively, the algebra of creation
operators on it.
\item[${\C[H]}$]The complex Hilbert space~$H$ is the Hilbert space of
single-particle states for a charged particle, and the algebra~$\C[H]$
is dense in the subspace of Fock space not including any
antiparticles. As an algebra,~$\C[H]$ is the algebra of creation
operators of particles, but it contains no creation operators of
antiparticles.
\item[${\C[H\oplus H^*]}$]The dual space~$H^*$ is the Hilbert space of
the single-antiparticle states associated to~$H$, and~$\C[H\oplus
H^*]$ is dense in the full Fock space of particles and
antiparticles. The isomorphism~$\C[H\oplus
H^*]\simeq\C[H_1^\sharp\oplus H_2^\sharp]$ is associated with the
possibility of representing a complex charged fields by a pair of real
neutral fields, and conversely. Again, as an algebra, this is the
algebra of creation operators of one species of charged particles and
antiparticles, or two species of truly neutral particles.
\item[${\C[H_1\oplus H_2]}$]This is dense in the Fock space of two
charged particles which are not antiparticles of each other, and it
does not include the antiparticle states. As an algebra, it is the
algebra of creation operators of two charged particles, with no
creation operators for antiparticles.
\end{description}

\section{The weak topology}\label{sec:WkTop}

In order to use tools from analysis it is necessary that all algebras
under consideration have a topology making all the algebra operations
continuous, and that the involution on a $*$-algebra be continuous,
too. This is not much of a restriction, since any vector space~$V$ can
be given the (locally convex) weak topology~$\omega(V,V^*)$ induced by
its algebraic dual~$V^*$, which then coincides with the topological
dual. Not only that, but every linear map~$f\colon V\to W$ is
continuous with respect to the weak topologies on~$V$ and~$W$. Indeed,
$f\colon V\to W$ is continuous with respect to the weak topology
on~$W$ if, and only if, $g\circ f\colon V\to\C$ is continuous for
all~$g\colon W\to\C$, but the weak topology on~$V$ makes every linear
functional on it, and in particular~$g\circ f$, continuous by
definition. Since the multiplication and unit maps are linear, the
weak topology provides a functor from algebras with homomorphisms to
locally convex algebras with continuous algebra homomorphisms. The
same applies to Abelian algebras and $*$-algebras, and the weak
topology defines functors from Abelian algebras or $*$-algebras to
locally convex Abelian algebras or $*$-algebras. Because the weak
topology can be put on every algebra in a way that makes all
homomorphisms continuous, we have the following commutative diagram of
functors:
$$
\xymatrix{\cat{AbAlg}\ar@<.5ex>[r]^{\mathrm{ F}}\ar[d]_{\omega(-,-^*)}&\cat{Ab{}Alg}^*\ar@<.5ex>[l]^{\mathrm{ U}}\ar[d]^{\omega(-,-^*)}\\
          \cat{LCAbAlg}\ar@<.5ex>[r]^{\mathrm{ F}}&\cat{LCAb{}^*Alg}\ar@<.5ex>[l]^{\mathrm{ U}}\\}
$$
On the bottom row of this diagram, the ``underlying''
functor~$U\colon{LCAb{}^*Alg}\to{LCAbAlg}$ takes any locally convex
$*$-algebra~$B$ to its underlying algebra~$U(B)$ with the weak
topology defined by the collection of all~$U(f)$, where~$f\colon
B\to\C$ is a $*$-homomorphism. Given a locally convex algebra~$A$, the
``free'' locally convex $*$-algebra~$F(A)$ is the free $*$-algebra
with the weak topology defined by all the~$*$-homomorphisms $F(f)$
where $f\colon A\to\C$ is a homomorphism.

\subsection{Examples}

Algebras of polynomials on finitely many variables, such as~$\C[z]$,
are isomorphic as vector spaces to the space of complex sequences with
finitely many nonzero entries, usually denoted~$c_{00}$. A linear
functional on this space assigns to any sequence a linear combination
of its (finitely many) nonzero entries, and the algebraic
dual~$c_{00}^*$ is isomorphic to the space of all sequences~$l_0$,
with unrestricted complex coefficients. The space~$c_{00}$ is weakly
complete, since for any~$a\in\overline{c_{00}}$ there is an~$\alpha\in
l_0$ obtained by replacing each nonzero element of~$a$ with its
inverse, so that~$\alpha(a)$ is the number of nonzero elements
of~$a$. Since~$a\in\overline{c_{00}}$ and~$\alpha\in c_{00}^*$,
$\alpha(a)$ must be finite, and so~$a\in c_{00}$ already.

In the case of~$\C[H]$, the set of generators of the polynomial
algebra is not just any infinite set, but it forms a Hilbert space. We
have a homogeneous decomposition $\C[H]\simeq\bigoplus_{n\ge
0}H^{\odot n}$, where~$H^{\odot n}$ (the symmetric tensor power
of~$H$) has a natural inner product derived from that of~$H$, and any
element of~$\C[H]$ is finite linear combination of monomials. We
define~$\C[H]^*\colon =\prod_{n\ge 0}(H^{\odot n})^*$, that is, each
element of~$\C[H]^*$ consists of one element of the topological
dual~$(H^{\odot n})^*$ for each~$n$. The resulting weak topology
makes~$\C[H]$ complete, as in the finitely-generated case.

\section{The Gel'fand spectrum}\label{sec:spectrum}

If~$A$ is an Abelian algebra or $*$-algebra~$A$ with the weak
topology, we define its Gel'fand spectrum, denoted~$\Delta_A$, as the
collection of morphisms from it into the complex numbers in the
appropriate category. Specifically, if~$A$ is a commutative algebra
with a topology, we define the Gel'fand spectrum of~$A$ to be the
collection of all continuous algebra homomorphisms into the complex
numbers. In symbols,
$$
\Delta_A=\cat{LCAbAlg}(A,\C).
$$
Similarly, if~$A$ is a commutative $*$-algebra with a topology, we
define its Gel'fand spectrum~$\Delta_A$ to be the collection of all
continuous $*$-algebra $*$-homomorphisms into the complex numbers, or
$$
\Delta_A=\cat{LCAbAlg^*}(A,\C).
$$
These are both instances of hom-sets so, by abstract nonsense, they
are contravariant functors to~$\cat{Set}$, meaning that algebra
homomorphisms induce natural set maps going between the spectra in the
opposite direction. Precisely, if~$f\colon A\to B$ is an continuous
homomorphism (or $*$-homomorphism) between Abelian algebras (or
$*$-algebras), then there is a function
$\Delta_f\colon\Delta_B\to\Delta_A$ given by~$\Delta_f(p)=p\circ
f\colon A\to\C\in\Delta_A$ for any continuous homomorphism (or
$*$-homomorphism)~$p\in\Delta_B$. In the literature, the Gel'fand
spectrum is normally defined as the collection of maximal ideals. The
functorial definition given here is much more restrictive, and it
coincides with the usual definitions only for~$C^*$-algebras or normed
algebras. We discuss this in greater detail in the next section.

Consider now the evaluation map
$$
\begin{array}{rccc}
e\colon&\Delta_A\times A&\to&\C\\ &(p,a)&\mapsto&p(a).
\end{array}
$$
Equivalent to this is the Gel'fand transform, which associates to each
element~$a\in A$ the function~$e(~,a)\colon\Delta_A\to\C$. The
Gel'fand transform 
$$
\begin{array}{rccc}
\hat{~}\colon &A&\to&\C^{\Delta_A}\\
&a&\mapsto&e(~,a)
\end{array}
$$
is an algebra homomorphism (or $*$-homomorphism) into the $*$-algebra
of all complex functions on~$\Delta_A$ (with pointwise complex
conjugation as involution).

The evaluation map induces a natural topology on~$\Delta_A$, namely
the weakest topology making every $\hat a\in\hat A$ continuous. Note
that~$\Delta_A$ is a subset of the dual~$A^*$, and that the spectral
topology just defined is the same as the one induced on~$\Delta_A$ as
a subset of~$A^*$ with the weak${}^*$ topology. Since the weak${}^*$
topology separates points---because given two different homomorphisms
$p,q\in\Delta_A$, there must be an $a\in A$ such that $\hat
a(p)=p(a)\neq q(a)=\hat a(q)$---,~$\Delta_A$ is Hausdorff.

\begin{proposition}[Completeness]
The spectrum~$\Delta_A$ is a weak${}^*$ closed subset of~$A^*$.
\end{proposition}

\paragraph[Completeness]{Proof}

This argument is essentially the first half of the proof of the
Banach--Alaoglu theorem~\cite[\S 3.15]{rudin91}.

Let~$\{p_\lambda\}_{\lambda\in\Lambda}$ be a net in~$\Delta_A$
converging in the weak${}^*$ topology; \ie, for every $a\in A$,
$p_\lambda(a)\to p(a)$ for some $p(a)\in\C$. Then,
$$
p(a+b)-p(a)-p(b)=\bigl[p(a+b)-p_\lambda(a+b)\bigr]+\bigl[p_\lambda(a)-p(a)\bigr]+\bigl[p_\lambda(b)-p(b)\bigr]
$$
implies that~$p$ is linear, so $p\in A^*$. Also,
\begin{eqnarray*}
p(ab)-p(a)p(b)&=&\bigl[p(ab)-p_\lambda(ab)\bigr]+\bigl[p_\lambda(a)-p(a)\bigr]\bigl[p_\lambda(b)-p(b)\bigr]+\\
&&+\bigl[p_\lambda(a)-p(a)\bigr]p(b)+p(a)\bigl[p_\lambda(b)-p(b)\bigr]
\end{eqnarray*}
together with the trivial observation that $p_\lambda(1)=1$ for
all~$\lambda\in\Lambda$ so $p(1)=1$, implies that~$p$ is an algebra
homomorphism. This completes the proof in the category of algebras and
continuous algebra homomorphisms. If~$A$ is a $*$-algebra,
$$
p(a^*)-\overline{p(a)}=\bigl[p(a^*)-p_\lambda(a^*)\bigr]+\bigl[\,\overline{p_\lambda(a)}-\overline{p(a)}\,\bigr]
$$
implies that~$p$ is a $*$-algebra $*$-homomorphism. \qed

The second part of the proof of the Banach--Alaoglu theorem~\cite[\S
3.15]{rudin91} provides a characterization of compact subsets of the
spectrum.

\begin{proposition}[Compactness]
\label{compact_spectrum}
With respect to the weak${}^*$ topology on~$\Delta_A$, a closed subset
$F\subseteq\Delta_A$ is compact if, and only if, every $\hat a\in\hat
A$ is bounded on it.
\end{proposition}

\paragraph[Compactness]{Proof}

\begin{itemize}
\item[$\Rightarrow$)] the continuous image of a compact set is
compact, and compact sets of~$\C$ are bounded; and
\item[$\Leftarrow$)] we can use~$\hat A$ to embed~$F$ homeomorphically
as a closed subset of a cube which is compact by Tychonoff's theorem,
and Hausdorff. \qed
\end{itemize}

It follows that a subset of the spectrum has compact closure if, and
only if, every~$\hat a\in\hat A$ is bounded on it; and a point of the
spectrum has a basis of compact neighbourhoods if, and only if, it has
an neighbourhood on which every~$\hat a\in\hat A$ is bounded. Also,
this result implies that the restriction of~$\hat A$ to a compact
subset of~$\Delta_A$ is a normed algebra.

We now turn to the question whether the set
map~$\Delta_f\colon\Delta_B\to\Delta_A$ defined above is a continuous
map with respect to the weak${}^*$ topologies on~$\Delta_A$
and~$\Delta_B$. This is all that is required to show that~$\Delta$ is
a functor not only into~$\cat{Set}$, but into~$\cat{Top}$.

\begin{proposition}[Functoriality]
If~$f\colon A\to B$ is a continuous homomorphism (or $*$-homomorphism)
of Abelian algebras (or $*$-algebras), then the set map
$$
\begin{array}{rccc}
\Delta_f\colon&\Delta_B&\to&\Delta_A\\
&p&\mapsto&fp\\
\end{array}
$$
is continuous with respect to the weak${}^*$ topologies on~$\Delta_A$
and~$\Delta_B$.
\end{proposition}

\paragraph[Functoriality]{Proof}

The weak${}^*$ topology on~$\Delta_A$ admits a sub-base consisting of
sets of the form~$U=\hat a^{-1}(G)$, where~$a\in A$ and~$G$ is open
in~$\C$. We need to show that~$V=(\Delta_f)^{-1}(U)$ is open with
respect to the weak${}^*$ topology on~$\Delta_B$. In fact, a stronger
statement is true, namely, $(\Delta_f)^{-1}(U)=\hat b^{-1}(G)$ where
$b=f(a)$. Indeed,~$p\colon B\to\C$ is in~$V$ if, and only
if,~$\Delta_f(p)=fp\in U$, that is, $p\bigl(f(a)\bigr)\in G$ or,
equivalently,~$\hat b(p)\in G$. \qed

\begin{proposition}[Separation and regularity]
With respect to the weak${}^*$ topology, $\Delta_A$ is a Tychonoff
(completely regular, Hausdorff) space. 
\end{proposition} 

\paragraph[Separation and regularity]{Proof}

The topology on~$\Delta_A$ is the weak topology defined by the complex
functions~$\hat a\in\hat A$. However, the same topology is obtained if
the image is considered to be not the complex plane, but the complex
sphere, which is compact metric and so Tychonoff. By means of the
family of all~$\hat a\in A$,~$\Delta_A$ can be homeomorphically
embedded as a subset of a product of Tychonoff spaces, and so is a
Tychonoff space~\cite[\S 14]{willard}. \qed

The situation is this:
$$
\xymatrix{\cat{LCAbAlg}\ar@<.5ex>[r]^{ F}\ar[d]_{\Delta_\cdot}&\cat{LCAb{}^*Alg}\ar[dl]^{\Delta_\cdot}\ar@<.5ex>[l]^{ U}\\
          \cat{Tych}\\}
$$
The diagram commutes in one direction only, namely, for any locally
convex Abelian algebra~$A$, it is true that~$\Delta_{ F(A)}=\Delta_A$
because each continuous algebra homomorphism~$f\colon A\to B$ extends
to a unique continuous $*$-algebra $*$-homomorphism~$F(f)\colon
F(A)\to F(B)$ whose restriction to~$A$ is precisely~$f$. On the other
hand, if~$A$ is a general locally convex Abelian $*$-algebra,
$\Delta_A\not\simeq\Delta_{ U(A)}$. However, the next best thing is
true: there is a natural
transformation~$j\colon\Delta_\cdot\Rightarrow\Delta_{ U(\cdot)}$
associated to the fact that every continuous $*$-algebra
$*$-homomorphism is an ordinary continuous homomorphism of the
underlying algebra. In other words,

\begin{proposition}[Naturality]
If~$A$ is any locally convex $*$-algebra and~$ U(A)$ is the underlying
locally convex algebra, there is a continuous inclusion
map~$j_A\colon\Delta_A\to\Delta_{ U(A)}$ such that, for every
continuous $*$-algebra $*$-homomorphism $f\colon A\to B$ the following
diagram commutes
$$
\xymatrix{A\ar[d]^f & \Delta_A\ar[r]^{j_A} & \Delta_{ U(A)}\\
          B & \Delta_B\ar[r]^{j_B}\ar[u]^{\Delta_f} &
          \Delta_{ U(B)}\ar[u]^{\Delta_{ U(f)}}}
$$
\end{proposition}

\paragraph{Proof}

Recall that $ U(f)\colon U(A)\to U(B)$ is the continuous algebra
homomorphism between the underlying locally convex algebras associated
to~$f\colon A\to B$. Recall also that~$\Delta_f\colon
\Delta_B\to\Delta_A$ is the continuous map obtained by composing
with~$f$; that is, if~$p\colon B\to\C$ is in~$\Delta_B$, then
$\Delta_f(p)=fp\colon A\to\C$ is in~$\Delta_A$. Similarly, $\Delta_{
U(f)}\colon \Delta_{ U(B)}\to\Delta_{ U(A)}$ is the continuous
function mapping~$p\colon U(B)\to\C$ to $ U(f)p\colon U(A)\to\C$.

To show that the diagram commutes, let~$p\colon B\to\C$ be a
continuous $*$-algebra $*$-homomorphism. Then, $\Delta_f(p)=fp\colon
A\to\C$, and $j_A(fp)\colon U(A)\to\C$ is the associated continuous
algebra homomorphism. On the other hand, $j_B(p)\colon U(B)\to\C$ is
the continuous algebra homomorphism associated to~$p$, and~$\Delta_{
U(f)}\bigl(j_B(p)\bigr)= U(f)j_B(p)\colon U(A)\to\C$. It remains only
to show that $j_A(fp)= U(f)j_B(p)\colon U(A)\to\C$, but this is
because, as set maps, $j_A(fp)=fp$, $j_B(p)=p$, and~$ U(f)=f$. \qed

We can sum up the content of this section in the following theorem.

\begin{theorem}[Gel'fand spectrum]\label{thm:spectrum}
Let~$A$ be a locally convex algebra or $*$-algebra, and let its
Gel'fand spectrum~$\Delta_A$ be the hom-set~$\hom(A,\C)$ in the
appropriate category. Then,~$\Delta_A$ is a weak${}^*$-closed subset of
the topological dual~$A^*$ and inherits a Tychonoff space
topology. \qed
\end{theorem}

\subsection{Examples}

Consider the $*$-algebra~$A=\C[x]$ where~$x^*=x$. A $*$-algebra
$*$-homomorphism $p\colon A\to\C$ is uniquely determined by~$p(x)$,
which must be a real number since $p(x)=p(x^*)=\overline{p(x)}$. Other
than that,~$p(x)\in\R$ is unrestricted, and so~$\Delta_A\simeq\R$.

Similarly, it can be shown that the spectrum of~$B=\C[z]$
is~$\Delta_B\simeq\C$. Since~$\R$ is strictly contained in~$\C$
and~$B= U(A)$, we have an example of how~$\Delta_{U
(A)}\neq\Delta_A$.

Next we consider the spectrum of~$C=\C[z,z^*]$. An algebra
$*$-homo\-morph\-ism~$p\colon C\to\C$ is determined by~$p(z)\in\C$, and
the condition that~$p(z^*)=\overline{p(z)}$ does not restrict the
possible value of~$p(z)$. Hence,~$\Delta_C\simeq\C$. This is expected,
as~$C= F(B)$ and we know that~$\Delta_{ F(B)}=\Delta_B$.

We can use these three examples to illustrate a principle: real
analysis is all about $*$-algebras, and complex analysis is all about
algebras. Also, analysis on the complex plane is done by going back
and forth between the algebra~$\C[z]$ and the $*$-algebra~$\C[z,z^*]$,
whose spectra are both isomorphic to~$\C$. The difference is
that~$\C[z,z^*]$ is used to study the structure of~$\C$ as a
two-dimensional real manifold, while~$\C[z]$ is used to study the
structure of~$\C$ as a one-dimensional complex manifold. In complex
analysis, nominally one is studying holomorphic functions, which are
limits of polynomials in~$\C[z]$. However, often one needs to use the
real and imaginary parts, which live in~$\C[z,z^*]\simeq\C[x,y]$. A
case in point is the Cauchy--Riemann equations~$\partial
f/\partial\bar z=0$, which characterizes the image of~$\C[z]$ inside~$
F\bigl(\C[z]\bigr)=\C[z,\bar z]$. In other words, the following sequence
is exact:
$$
\xymatrix{\C[z]\ar[r]^{F}&\C[z,z^*]\ar[r]^{\partial/\partial
z^*}&\C[z,z^*]} 
$$

The case of~$\C[H]$, where~$H$ is a Hilbert space, is interesting
because its spectrum is not locally compact. Indeed, just as in the
case of polynomials on finitely many generators, an algebra
homomorphism~$p\colon\C[H]\to\C$ is uniquely determined
by~$\left.p\right|_H\in H^*$, and so~$\Delta_{\C[H]}\simeq H^*$, with
the weak${}^*$ topology. We know that locally compact, Hausdorff
topological vector spaces must be finite-dimensional, so in this case
the spectrum is not locally compact. Incidentally, since every unital
Banach algebra has compact spectrum, this shows that the algebra of
creation operators on Fock space cannot be a Banach algebra.

\section{The Gel'fand transform}\label{sec:trans}

We now study in detail the Gel'fand transform, which is the algebra
homomorphism (or $*$-algebra $*$-homomorphism) given by
$$
\begin{array}{rccc}
\hat{~}\colon &A&\to&C(\Delta_A)\\
&a&\mapsto&e(~,a)
\end{array}
\qquad\hbox{such that}\quad
\begin{array}{rccc}
\hat a=e(~,a)\colon &\Delta_A&\to&\C\\
&p&\mapsto&p(a)
\end{array}
$$
where~$C(\Delta_A)\subseteq\C^{\Delta_A}$ denotes the $*$-algebra of
continuous complex functions on~$\Delta_A$ or its underlying
algebra. To fully understand this homomorphism we need to characterize
its kernel and its image.

\subsection{Ideals and homomorphisms}

The kernel of a homomorphism of Abelian algebras is an ideal, that is,
closed under addition and preserved by multiplication by elements of
the algebra. Conversely, the quotient of an Abelian algebra by an
ideal is an Abelian algebra homomorphism. The kernel of a
$*$-homomorphism is closed under the involution and, if an ideal is
closed under the involution the quotient is a $*$-homomorphism.

All nilpotent elements of~$A$ must be in the kernel of the Gel'fand
transform, as the equation~$a^n=0$ translates into the complex
equation~$p(a)^n=0$, for all~$p\in\Delta_A$, which is equivalent
to~$p(a)=0$ for all~$p\in \Delta_A$, or~$\hat a=0$. Although the
nilpotent elements form an ideal, it is possible that the kernel of
the Gel'fand transform contains other elements. If the Gel'fand
transform is one-to-one, we say the algebra~$A$ is semisimple.

We have defined the Gel'fand spectrum as the collection of continuous
algebra homomorphisms (or $*$-algebra $*$-homomorphisms) into~$\C$. We
call the kernels of these homomorphisms Gel'fand ideals, and they are
characterized by being closed, codimension-$1$ ideals and, in the case
of~$*$-algebras, closed under the involution. The kernel of the
Gel'fand transform, called the Gel'fand radical of~$A$, consist of
precisely those~$a\in A$ on which every~$p\in\Delta_A$ vanishes. Being
the intersection of all the Gel'fand ideals, it is a closed ideal and,
if~$A$ is a $*$-algebra, it is closed under the involution. In sum,
\begin{description}
\item[Gel'fand ideal] An ideal~$I$ in~$A$ of codimension~$1$, closed
if~$A$ has a topology and closed under the involution if~$A$ is a
$*$-algebra. The quotient~$A/I$ is~$\C$, and the quotient map is
a~$*$-homomorphism if~$A$ is a~$*$-algebra.
\item[Gel'fand radical] The Gel'fand radical~$R$ is the kernel of the
Gel'fand transform, and it is the intersection of all Gel'fand
ideals. If~$A$ has a topology, the Gel'fand radical is closed; if it
has an involution, the Gel'fand radical is closed under it.
\end{description}

Let us now analyze in more detail the difference between maximal
ideals and Gel'fand ideals or, equivalently, the difference between
the Gel'fand radical and the Jacobson radical. Recall the following
concepts from commutative algebra~\cite{atiyah69}:
\begin{description}
\item[Maximal ideal] A proper ideal~$I$ in~$A$ is maximal iff it is
maximal among proper ideals with respect to inclusion. The
algebra~$A/I$ is a field. Every Gel'fand ideal is maximal, but maximal
ideals may fail to be closed or have codimension~$1$.
\item[Prime ideal] A proper ideal~$I$ in~$A$ is prime iff $a,b\not\in
I$ implies $ab\not\in I$. In the algebra~$A/I$, the product of nonzero elements
is nonzero. Every maximal ideal is prime.
\item[Radical ideal] If~$I$ is an ideal in~$A$, the radical of~$I$ is
the ideal
$$
\rad(I)=\{a\in A\mid\exists n>0,a^n\in I\}.
$$
It is the intersection of the prime ideals containing~$I$.
\item[Jacobson radical] The Jacobson radical of~$A$ is the ideal~$J$
obtained by taking the intersection of all maximal ideals of~$A$, and
it is contained in the Gel'fand radical.  
\item[Nilradical] The nilradical of~$A$ is the ideal~$N$ consisting of
all nilpotent elements of~$A$ (i.e., the radical ideal of the zero
ideal). It is the intersection of all prime ideals. The algebra~$A/N$
has no nilpotent elements. It is contained in the Jacobson radical.
\end{description}

We have already indicated that the Gel'fand spectrum is usually
defined as the collection of all maximal ideals, the implication being
that the Gel'fand radical coincides with the Jacobson radical. This is
because maximal ideals of a normed algebra are closed, and
because~$\C$ is the only normed field extension of~$\C$ (the
Gel'fand--Mazur theorem), every maximal ideal of a normed algebra is a
Gel'fand ideal. In the presence of an involution, the construction
only works for~$C^*$-algebras because only then it can be proved that
every homomorphism is a~$*$-homomorphism.

Since every $*$-homomorphism is a homomorphism, for $*$-algebras what
we have called the Gel'fand spectrum is in general smaller than
usually defined, and the Gel'fand radical larger than usual. Sometimes
the spectrum of an algebra is defined as the collection of all algebra
homomorphisms into~$\C$ (also called characters), continuous or not,
irrespective of whether the algebra under consideration has an
involution. Because of the inclusion of ordinary homomorphisms in the
spectrum of a $*$-algebra, it can happen that the Gel'fand transform
is not a~$*$-homomorphism. This is fixed by removing, as we do, from
the spectrum of a $*$-algebra all the homomorphisms which are not
$*$-homomorphisms.

In ring theory it is remarked that the spectrum of prime ideals is
functorial because the inverse image of a prime ideal by a
homomorphism is a prime ideal, but not so for maximal ideals so the
maximal spectrum is not functorial. The point of our redefinition of
the Gel'fand transform is to show that, by insisting on a functorial
definition that extends beyond the realm of $C^*$-algebras, some of
the important conclusions of the Gel'fand--Na\u\i{}mark theorem can
also be extended.

\subsection{Topologies on~$C(\Delta)$}

The space~$C(\Delta_A)$ has a weak topology making all the evaluation
maps continuous, which can be easily seen to be associated to
pointwise convergence on the spectrum; the image of the Gel'fand
transform~$\hat A$ inherits this topology. On the other hand, since
the Gel'fand radical~$G$ is a closed ideal of the locally convex
algebra~$A$, the image of the Gel'fand transform~$\hat A\simeq A/G$
has a locally convex quotient topology. These two topologies on~$\hat
A$ coincide. Note, however, that the space of continuous functions is
rarely closed under the topology of pointwise convergence. A stronger
topology is needed to make the algebra~$C(\Delta_A)$ closed, but then
it is no longer obvious that the Gel'fand transform is continuous.

The natural stronger topology on~$C(\Delta_A)$ is the compact-open
topology (\ie, uniform convergence on compact sets), which is the
locally convex topology defined by the seminorms
$$
|f|_K=\sup_{p\in K}|f(p)|,
$$
where~$K$ is any compact subset of~$\Delta_A$, and the algebra
operations are continuous with respect to this topology. Since~$\hat A$
is a subalgebra of~$C(\Delta_A)$, it inherits the compact-open
topology. The original (weak) topology on~$A$ is strictly weaker than
the compact-open topology on~$\hat A$ unless the only compact subsets
of the spectrum~$\Delta_A$ are the finite subsets.

The compact-open topology is natural in another, more interesting
sense, and that is the existence of a Stone--Weierstrass theorem for
Tychonoff spaces (see~\cite[\S44B]{willard} for a sketch of the
proof). Since the spectrum~$\Delta_A$ of any algebra~$A$ is a
Tychonoff space, it follows that~$\hat A$ is dense in~$C(\Delta_A)$
with the compact-open topology.

\begin{proposition}[Stone--Weierstrass for Tychonoff spaces]
If~$\Delta$ is a Tychonoff space and~$A$ is a~$*$-subalgebra
of~$C(\Delta)$ which separates points of~$\Delta$, and contains the
constant functions, then~$A$ is dense in~$C(\Delta)$ with the
compact-open topology. \qed
\end{proposition}

The conclusion is that the following diagram of functors commutes in
both directions.
$$
\xymatrix{&\cat{LCAb{}^*Alg}\ar[dl]_{\Delta_\cdot}\ar[d]^{A\mapsto\bar{\hat
          A}}\\
          \cat{Tych}\ar@<.5
          ex>[r]^{C(-)}&\cat{AbLC^*Alg}\ar@<.5
          ex>[l]^{\Delta_\cdot}\\}
$$
That is, because of the Stone--Weierstrass theorem for Tychonoff
spaces, the closure of~$\hat A$ is the space of continuous functions on
the Gel'fand spectrum; and then there is the rather trivial
observation that the spectrum of~$C(\Delta)$ is
precisely~$\Delta$. This last observation implies that the Gel'fand
functor from Abelian $C^*$-algebras to compact Hausdorff spaces, and
its inverse, extend to functors between Abelian $LC^*$-algebras and
Tychonoff spaces:
$$
\xymatrix{\cat{AbC^*Alg}\ar@<.5ex>[r]^{\Delta}\ar[d]&\cat{CompT_2}\ar@<.5ex>[l]^{C(\cdot)}\ar[d]\\
          \cat{AbLC^*Alg}\ar@<.5ex>[r]^{\Delta}&\cat{Tych}\ar@<.5ex>[l]^{C(\cdot)}}
$$

We can summarize the content of this section in the following
theorem. For lack of a better name, we call the algebra of continuous
complex functions on a Tychonoff space an~``$LC^*$-algebra'', for
``locally convex'' and ``locally~$C^*$''.

\begin{theorem}[Generalized Gel'fand--Na\u\i{}mark
theorem]\label{thm:transform}  
If~$\Delta_A$ is the Gel'fand spectrum of a semisimple, locally convex
$*$-algebra~$A$, the Gel'fand transform is $*$-homomorphism of~$A$
into a dense $*$-subalgebra of~$C(\Delta_A)$, the~$*$-algebra of
continuous complex functions with the compact-open topology.
\end{theorem}

\subsection{Examples}

In the case of the $*$-algebra~$A=\C[x]$, the results of this section
translate into the fact that complex polynomials on~$\R$ (the Gel'fand
spectrum of~$A$) are dense in the space of all continuous functions
from~$\R$ to~$\C$ with the compact-open topology. Similarly, the space
of polynomials~$\C[z,z^*]$ is dense in the continuous functions on~$\C$
with the compact-open topology.

In the infinite-dimensional case, we get the more interesting result
that~$\C[H,H^*]$ is dense (with the compact-open topology) in the
space of all continuous complex functions on the Hilbert
space~$H$. This goes a long way towards reducing nonlinear analysis on
Hilbert spaces to algebra.

\section{The Gel'fand--Na\u\i{}mark--Segal construction}\label{sec:GNS}

The Gel'fand--Na\u\i{}mark--Segal theorem is based on the concept of a
state on a~$C^*$-algebra, which in the commutative setting has the
interpretation of a classical expectation value on a family of bounded
random variables. Since the definition of state does not require the
algebra to be a~$C^*$-algebra, it applies without modification to our
setting. An intuitively appealing characterization of states which
uses the Riesz representation theorem is that any state on
a~$C^*$-algebra is realized as a Borel probability measure on the
Gel'fand spectrum. As we have seen, all that is lost when dropping
the~$C^*$ hypothesis is the compactness of the spectrum, but the next
best result is true: the Gel'fand spectrum, if correctly defined, is
always a Tychonoff space.

A state~$E$ on a $*$-algebra~$A$ is a positive, normalized,
compact-open continuous linear functional on~$\hat A$. That is:
\begin{itemize}
\item $E\in A^*$,
\item $E(1)=1$, and
\item $E(\hat a^*\hat a)\ge 0$ for all $a\in A$. 
\item there are compact subsets $K_1,\ldots,K_n\subseteq\Delta_A$ and 
      positive numbers $C_1,\ldots,C_n$ such that
      $\max\bigl\{C_i|\hat a|_{K_i}\bigr\}<1$ implies
      $\bigl|E(\hat a)\bigr|<1$;
\end{itemize}
The compact-open continuity condition is
equivalent to
\begin{itemize}
\item there is a compact subsets $K\subseteq\Delta_A$ and a
      positive number $C$ such that
      $C|\hat a|_K<1$ implies
      $\bigl|E(\hat a)\bigr|<1$
\end{itemize}
(just let $C=\max\{C_i\}$ and $K=\cup\{K_i\}$). By the Riesz
representation theorem, any such linear functional is of the form
$$
E(a)=\int_K \hat a(p)\mu_E(dp)
$$
where $\mu_E$ is a Borel probability measure on $K$. In other words, a
state on~$A$ is a state on the $C^*$-algebra $\left.A\right|_K$. In
this way, once the Gel'fand spectrum of an algebra or $*$-algebra is
known, an ample supply of states becomes available. 

Given a state~$E$ on the $*$-algebra~$A$, one can define an inner
product on~$A$ by the formula $\langle a,b\rangle\colon =E(a^*b)$ for
all $a,b\in A$. To obtain a Hilbert space one must complete~$A$ with
respect to the inner product, and take the quotient by the ideal of
zero-norm elements of~$A$. This is the Gel'fand--Na\u\i{}mark--Segal
construction. As in the case of $C^*$-algebras~\cite[\S
12.41]{rudin91}, given that for every~$\hat a\in \hat A$ it is always
possible to find a state that does not vanish on it, one can find a
(possibly non-separable) Hilbert space on which~$\hat A$
and~$\C(\Delta_A)$ are faithfully represented as algebras of unbounded
operators.

If $\Delta_A$ is not compact, Borel probability measures on $\Delta_A$
are associated with densely defined states, meaning positive,
normalized, linear functionals on~$C(\Delta_A)$ which are densely
defined and not necessarily continuous in the compact-open
topology. The GNS construction can be carried out normally in that
case, with due attention being paid to subtleties about domains of
unbounded operators.

\subsection{Examples}

It might be surprising that states on the algebra of polynomials must
be compactly-supported measures on the spectrum. Where did the
ubiquitous Gaussian measure go? The answer is that the Gaussian
probability density function has an inverse which can be approximated
by polynomials uniformly on compact sets, and that means that the
Gaussian measure cannot be a continuous linear functional with respect
to the compact-open topology. However, the integral of the Gaussian
density times any polynomial is finite, so the Gaussian measure is a
densely-defined state. A similar argument holds for the
two-dimensional Gaussian and the algebra~$\C[z,z^*]$. With due care,
the Gaussian measure and other measures with non-compact support are
no harder to deal with than compactly-supported measures.

\section{Conclusions}

In this paper we have shown that the Gel'fand--Na\u\i{}mark theorem
generalizes from~$C^*$-algebras to any semisimple~$*$-algebra, and
that the Gel'fand spectrum and Gel'fand transform are well-behaved for
virtually any algebra or~$*$-algebra. The key to obtaining these
results is to define the Gel'fand spectrum in a manifestly functorial
way which, if nothing else, shows the power of elementary category
theory as an aid to generalization and to the formulation of the right
definitions. 

The ability to generalize the Gel'fand--Na\u\i{}mark theorem to
essentially arbitrary algebras is important for applications in
physics and probability theory, where often the requirement that
observables be bounded seems rather unnatural and can be justified
only on the grounds of mathematical convenience. It would be desirable
to more fully illustrate the usefulness of our results with problems
where the~$C^*$-algebraic formulation of probability theory is awkward
because of the essential presence of unbounded random variables. On
the mathematical side, the question remains open whether there is an
intrinsic characterization of semisimple, locally convex~$*$-algebras.

\section{Acknowledgements}

I would like to thank my advisor, John C. Baez, for encouragement in
the writing of this paper. I am also indebted to fellow student Toby
Bartels for help with category theory, and to Daniel Grubb for
comments posted on the USENET newsgroup sci.physics.research. Finally
I would like to recognize the Perimeter Institute for Theoretical
Physics, where this research was initiated, and Fotini Markopoulou for
inviting me to visit there; and the Department of Mathematics at the
University of California at Riverside for supporting my graduate
study.

\bibliography{algebras}
\bibliographystyle{alpha}

\end{document}